# POSITIVE RECURRENCE OF PROCESSES ASSOCIATED TO CRYSTAL GROWTH MODELS

By E. D. Andjel, M. V. Menshikov[1] and V. V. Sisko[2]

*Université de Provence, University of Durham and Universidade de São Paulo*

We show that certain Markov jump processes associated to crystal growth models are positive recurrent when the parameters satisfy a rather natural condition.

**1. Introduction.** Gates and Westcott studied some Markov processes representing crystal growth models. In these models particles accumulate on a finite set of sites. The first two theorems in the present paper study the model obtained when this set of sites is one-dimensional. To define the process associated to that model, let

$$r(a,b,c) = \begin{cases} \beta_2, & \text{if } b < \min\{a,c\}, \\ \beta_1, & \text{if } \min\{a,c\} \leq b < \max\{a,c\}, \\ \beta_0, & \text{if } \max\{a,c\} \leq b, \end{cases}$$

and for each $n \geq 2$, let $X^n(t) = (X_1^n(t), \ldots, X_n^n(t))$ be a Markov jump process on $\{\mathbb{Z}_+\}^n$ such that

$$X_i^n(t) \to X_i^n(t) + 1, \qquad 1 \leq i \leq n,$$

at rate

$$r(X_{i-1}^n(t), X_i^n(t), X_{i+1}^n(t)).$$

While these rates are well defined if $1 < i < n$, their expression for $i = 1$ or $i = n$ depends on which boundary conditions we adopt. In this paper we will consider two boundary conditions:

Received April 2005; revised November 2005.
[1]Supported in part by FAPESP (2004/13610–2).
[2]Supported by FAPESP (2003/00847–1).

(a) Zero boundary conditions: $X_0^n(t) = X_{n+1}^n(t) = 0$ for all $t \geq 0$.

(b) Periodic boundary conditions: $X_0^n(t) = X_n^n(t)$ and $X_{n+1}^n(t) = X_1^n(t)$ for all $t \geq 0$.

We now define for $1 \leq i \leq n$

$$\Delta_i X^n(t) = X_i^n(t) - X_{i+1}^n(t).$$

Let

$$\tilde{r}(u,v) = \begin{cases} \beta_2, & \text{if } \min\{u,v\} > 0, \\ \beta_1, & \text{if } \min\{u,v\} \leq 0 < \max\{u,v\}, \\ \beta_0, & \text{if } \max\{u,v\} \leq 0. \end{cases}$$

Then, $r(a,b,c) = \tilde{r}(a-b, c-b)$. Therefore, the process

$$(\Delta_1 X^n(t), \ldots, \Delta_n X^n(t))$$

is Markovian too.

Gates and Westcott in [4] considered this process under periodic boundary conditions and parameters $\beta$ such that

(1.1) $$0 < \beta_0 < \beta_1 < \beta_2.$$

There, they proved that it is positive recurrent if either $n \leq 4$ or

$$(n-1)^2 \beta_0 < \beta_2.$$

In an earlier paper [3], the same authors studied the case in which the parameters $\beta$ satisfy (1.1) and $\beta_1 = \frac{1}{2}(\beta_0 + \beta_2)$. For this case, they gave for any $n$ an explicit expression for a stationary measure of the process, thus, proving its positive recurrence.

In this paper we show that (1.1) is a sufficient condition for positive recurrence for all values of $n$. More explicitly, we prove the following:

THEOREM 1.1. *Let $\beta_0, \beta_1$ and $\beta_2$ satisfy (1.1). Then, for periodic boundary conditions and for any $n \geq 2$, the process $(\Delta_1 X^n(t), \ldots, \Delta_n X^n(t))$ is positive recurrent and its unique invariant measure has exponentially decaying tails.*

THEOREM 1.2. *Let $\beta_0, \beta_1$ and $\beta_2$ satisfy (1.1). Then, for zero boundary conditions and for any $n \geq 2$, the process $(\Delta_1 X^n(t), \ldots, \Delta_{n-1} X^n(t))$ is positive recurrent and its unique invariant measure has exponentially decaying tails.*

In the statement of the last theorem the process considered is

$$(\Delta_1 X^n(t), \ldots, \Delta_{n-1} X^n(t))$$

rather than

$$(\Delta_1 X^n(t), \ldots, \Delta_n X^n(t)),$$

because for zero boundary conditions, we clearly have $\lim_t \Delta_n X^n(t) = \infty$. Note that the former process is Markovian too.

In Section 3 we prove that all the processes considered and starting from $(0, \ldots, 0)$ are such that

$$\sup_{\substack{i=1,\ldots,n-1,\\ t\geq 0}} P(|\Delta_i X^n(t)| \geq k)$$

decays exponentially in $k$. It is then easy to conclude that Theorems 1.1 and 1.2 hold. The proof is first given for zero boundary conditions and relies on an induction on $n$, some coupling arguments and Doléans exponential martingales.

When $\beta_1 = \frac{1}{2}(\beta_0 + \beta_2)$, the invariant measure given in [3] is such that the distribution of $\Delta_i X^n$ is the same for all $n$. It is therefore natural to ask if, under condition (1.1), there is a lower bound of the constant associated with the exponential decay of the invariant measure that is valid for all $n$. Unfortunately, our proof does not provide such a bound and its existence remains an open problem.

In Section 2 we prove the positive recurrence of a class of discrete time Markov chains. The precise result is stated in Theorem 1.3. A consequence of this theorem is that the processes considered by Theorems 1.1 and 1.2 are positive recurrent (see Examples 1 and 2 below). This proof does not give any information about the decay of the invariant measure, but we believe that many readers will be interested in it because it is much shorter than the proof given in Section 3 and because it treats many models which are not covered by Theorems 1.1 and 1.2 (see Example 3).

To describe the class of Markov chains considered in Section 2, let $G = (V, E)$ be a connected nonoriented finite graph, where $V = \{1, \ldots, n\}$. The elements of $V$ are seen as different columns, on which particles accumulate with time. A Markov chain $\zeta_t, t \geq 0$ from that class will have as state space $Y = \mathbb{Z}^{n-1}$ and its transition probability matrix will by denoted by $Q$. The $i$th coordinate of the chain represents the difference between the number of particles in columns $i$ and $n$.

In Theorem 1.3 we impose some conditions on $Q$, using the following functions:

$$f_{ij}(y) = \begin{cases} y(i) - y(j), & \text{if } i < n, j < n, \\ y(i), & \text{if } i < n, j = n, \\ -y(j), & \text{if } i = n, j < n, \end{cases}$$

where $y \in Y$ and $1 \leq i, j \leq n$. Hence, in all cases $f_{ij}(y)$ is the algebraic difference between the number of particles in columns $i$ and $j$.

For $i = 1, \ldots, n$, let $e_i \in Y$ be the vector that we have to add to an element of $Y$ when column $i$ increases by one unit, that is,

$$e_1 = (1, 0, \ldots, 0), \ldots, e_{n-1} = (0, 0, \ldots, 1), e_n = (-1, \ldots, -1, -1).$$

We now state our third theorem:

THEOREM 1.3. *Consider a graph $G = (V, E)$ and a Markov chain $\zeta_t$ with state space $Y = \mathbb{Z}^{n-1}$. Let $\delta$ and $M$ be constants such that $\delta, M \in (0, 1)$. Suppose that, for the transition probability matrix $Q$ of $\zeta_t$, the following conditions hold:*

(i) $Q(x, y) = 0$ *unless* $y = x$ *or* $y = x + e_i$ *for some* $1 \leq i \leq n$.
(ii) $\inf_{y \in Y} Q(y, y + e_i) \geq \delta$ *for all* $1 \leq i \leq n$.
(iii) *Suppose that* $y \in Y$, $\{i, j\} \in E$, *and* $f_{ij}(y) > 0$. *Then*

$$Q(y, y + e_i) \leq Q(y, y + e_j).$$

(iv) *Suppose that* $y \in Y$, $i \in V$ *and for all $\ell$ such that $\{i, \ell\} \in E$, we have $f_{i\ell}(y) > 0$. Then*

$$Q(y, y + e_i) \leq Q(y, y + e_\ell) - M$$

*for any $\ell$ as above.*

*Then, the Markov chain $\zeta_t$ is positive recurrent.*

We end this section with three examples of Markov chains to which this theorem can be applied. All these examples can be seen as crystal growth models. The first and second examples treat the embedded chains appearing in Theorems 1.2 and 1.1, respectively. The third example does not require the graph to have a one-dimensional structure.

EXAMPLE 1. Adopting the same notation as in Theorem 1.2, consider the process

$$Z(t) = (Z_1(t), \ldots, Z_{n-1}(t)),$$

where

$$Z_i(t) = \sum_{k=1}^{n-i} \Delta_{n-k} X^n(t), \qquad i = 1, \ldots, n-1.$$

This process is Markovian because it is obtained applying a one to one map to another Markovian process. The embedded Markov chain of the process $Z(t)$ has as state space $Y = \mathbb{Z}^{n-1}$ and its probability transition matrix $Q$ is such that, for all $y \in Y$, we have $Q(y, z) = 0$ unless $z = y + e_i$ for some $i = 1, \ldots, n$.

Now, let $G = (V, E)$ be given by
$$V = \{1, \ldots, n\}, \qquad E = \{\{i, i+1\} : i = 1, \ldots, n-1\},$$
then, using (1.1), it is easy to check that this chain satisfies all the hypothesis of Theorem 1.3. Hence, it is positive recurrent. This implies that the process $Z(t)$ is positive recurrent because its rates are bounded. Since this process is obtained applying a one to one map to the process
$$(\Delta_1 X^n(t), \ldots, \Delta_{n-1} X^n(t)),$$
this last process is positive recurrent too.

EXAMPLE 2. In this example $(\Delta_1 X^n(t), \ldots, \Delta_n X^n(t))$ is as in Theorem 1.1. We start noting that in this case $\sum_{i=1}^n \Delta_i X^n(t) \equiv 0$. Therefore, $(\Delta_1 X^n(t), \ldots, \Delta_{n-1} X^n(t))$ is also Markovian and it suffices to show that this last process is positive recurrent. This is done almost exactly as in the previous example, the only difference being that now
$$E = \{\{i, i+1\} : i = 1, \ldots, n-1\} \cup \{1, n\}.$$

EXAMPLE 3. Suppose that condition (1.1) is satisfied and that $\beta_2 \leq 1$. Then, let $G = (V, E)$ be any admissible graph, that is, it is nonoriented, finite and connected. As in the other examples, we take $V = \{1, \ldots, n\}$. Then, define the matrix of transition probabilities as follows: for $y \in Y$ and $i = 1, \ldots, n$, let
$$Q(y, y + e_i) = \begin{cases} \beta_0/n, & \text{if } f_{i\ell}(y) > 0 \text{ for all } \ell \text{ such that } \{i, \ell\} \in E, \\ \beta_2/n, & \text{if } f_{i\ell}(y) \leq 0 \text{ for all } \ell \text{ such that } \{i, \ell\} \in E, \\ \beta_1/n, & \text{otherwise}, \end{cases}$$
and let
$$Q(y, y) = 1 - \sum_{i=1}^n Q(y, y + e_i).$$

It is now easy to see that the conditions of Theorem 1.3 are satisfied with $M = (\beta_1 - \beta_0)/n$ and $\delta = \beta_0/n$.

**2. Positive recurrence.** To prove Theorem 1.3, we will apply the following variation of Foster's theorem:

THEOREM 2.1. *Let $\zeta_t$ be an irreducible Markov chain on a countable state space $Y$. Then, $\zeta_t$ is positive recurrent if there exist a positive function $f$ defined on $Y$, a bounded strictly positive integer-valued function $k$ also defined on $Y$ and a finite subset $A$ of $Y$ such that the following inequalities hold:*

(2.1) $\qquad \mathrm{E}(f(\zeta_{t+k(\zeta_t)}) - f(\zeta_t)|\zeta_t = y) \leq -1, \qquad y \notin A,$

(2.2) $\qquad \mathrm{E}(f(\zeta_{t+k(\zeta_t)})|\zeta_t = y) < \infty, \qquad y \in A.$

This theorem follows easily from Theorem 2.2.4 of [2]. Although in that reference the Markov chain is also assumed to be aperiodic, this extra condition is not needed if we are only interested in positive recurrence instead of ergodicity.

Throughout this section we adopt the same notation as in Theorem 1.3, assume its hypothesis and let

$$p = |E| \qquad \text{(i.e., the cardinality of } E\text{)}.$$

We start defining the functions $f$ and $k$ to which we will apply Theorem 2.1. Let

$$f(y) = \sum_{\{i,j\} \in E} f_{ij}^2(y).$$

We need the following lemmas.

LEMMA 2.1. *For all $\{i,j\} \in E$ and all $y \in Y$, we have*

$$\mathrm{E}(f_{ij}^2(\zeta_{t+1}) - f_{ij}^2(\zeta_t)|\zeta_t = y) \leq 1.$$

PROOF. The lemma is a consequence of the following observations:

(i) $|f_{ij}(\zeta_{t+1}) - f_{ij}(\zeta_t)| \leq 1$ almost surely,
(ii) from condition (iii) of Theorem 1.3, it follows that, for $\zeta_t$ such that $|f_{ij}(\zeta_t)| \geq 1$, we have

$$\mathrm{P}(|f_{ij}(\zeta_{t+1})| = |f_{ij}(\zeta_t)| + 1) \leq \mathrm{P}(|f_{ij}(\zeta_{t+1})| = |f_{ij}(\zeta_t)| - 1). \qquad \square$$

LEMMA 2.2. *Suppose $C \geq \frac{1+p}{2M}$ and let*

$$D = \{y \in Y : |f_{ij}(y)| > C, \{i,j\} \in E\}.$$

*Then inequality (2.1) holds for $k \equiv 1$ and any $y \in D$.*

In words, the set $D$ is the set of all states of the Markov chain $\zeta_t$ such that, for every $i, j \in \{1, \ldots, n\}$ such that $i$ and $j$ are neighbor columns, we have that the absolute value of the difference between heights of columns $i$ and $j$ is greater than $C$.

PROOF OF LEMMA 2.2. For a given $y \in D$, let $u$ be in the set of the columns that have the maximal number of particles and let $v$ be such that $\{u,v\} \in E$. We write

$$\mathrm{E}(f(\zeta_{t+1}) - f(\zeta_t)|\zeta_t = y)$$
$$= \mathrm{E}(f_{uv}^2(\zeta_{t+1}) - f_{uv}^2(\zeta_t)|\zeta_t = y)$$
$$+ \sum_{\substack{\{i,j\} \in E, \\ \{i,j\} \neq \{u,v\}}} \mathrm{E}(f_{ij}^2(\zeta_{t+1}) - f_{ij}^2(\zeta_t)|\zeta_t = y).$$

It is easy to see that $f_{uv}(y) = |f_{uv}(y)| > C$. From condition (iv) of Theorem 1.3, it follows that the first term in the right-hand side is bounded above by $-2CM + 1$. By Lemma 2.1, the second term is less than or equal to $p - 1$. Putting the bounds together, we complete the proof. $\square$

Note that the number of elements in the set $Y \setminus D$ is infinite. Therefore, we are not ready yet to apply Theorem 2.1 and finish the proof of Theorem 1.3.

Before the formulation of the last lemma we need some notation. Let

$$0 < C_1 < C_2 < \cdots < C_p < \infty.$$

The values of these constants will be determined later. Let

$$D_0 = \{y \in Y : |f_{ij}(y)| < C_p, \{i,j\} \in E\},$$
$$D_1 = \{y \in Y : |f_{ij}(y)| \geq C_1, \{i,j\} \in E\},$$
$$D_m = D'_m \cap D''_m, \qquad m = 2, \ldots, p,$$

where

$$D'_m = \{y \in Y : \text{ there exists } \{u,v\} \in E \text{ such that } |f_{uv}(y)| \geq C_m\}$$

and

$$D''_m = \{y \in Y : \text{ for any } \{i,j\} \in E \text{ we have } |f_{ij}(y)| \notin [C_{m-1}, C_m)\}.$$

Since $E$ has $p$ elements, we have

$$Y = D_0 \cup D_1 \cup D_2 \cup \cdots \cup D_p.$$

The following lemma constitutes the main step in the proof of Theorem 1.3.

LEMMA 2.3. *There exist positive constants $C_1, \ldots, C_p$ and positive integers $k_1, \ldots, k_p$ such that, for every $m \in \{1, \ldots, p\}$, inequality (2.1) holds for $y \in D_m$ if $k(y) = k_m$.*

PROOF. We take an increasing sequence of integers $C_m$ such that

$$C_1 \geq \frac{1+p}{2M},$$
$$C_m \geq \max\left\{pC_{m-1}, \frac{1 + p + p^2 C_{m-1}}{M \delta^{pC_{m-1}}}\right\}, \qquad m = 2, \ldots, p,$$

and let

$$k_m = \begin{cases} 1, & \text{if } m = 1, \\ 1 + pC_{m-1}, & \text{if } m = 2, \ldots, p. \end{cases}$$

Recall that the sets $D_m$ depend on the constants $C_m$. By Lemma 2.2, inequality (2.1) holds for $y \in D_1$ if $k(y) = k_1$.

Suppose that $y \in D_m$ for some $m = 2, \ldots, p$. Let $\ell$ be in the set of the columns that have the maximal number of particles. Consider the subgraph $\tilde{G} = (V, \tilde{E})$ of the graph $G = (V, E)$, where

$$\tilde{E} = \{\{i,j\} \in E : |f_{ij}(y)| < C_{m-1}\}.$$

By $O$ denote the largest connected subgraph of $\tilde{G}$ containing $\ell$.

Let us prove that $O \neq V$. Assume the converse. Recall that $|E| = p$. Then, for any $\{i,j\} \in E$, we have $|f_{ij}(y)| < pC_{m-1} \leq C_m$. This contradicts $y \in D'_m$.

Also, since the graph $G$ is connected, we see that there exists $\{u,v\} \in E$ such that $u \in O$, $v \notin O$, and $|f_{uv}(y)| \geq C_{m-1}$. Since $y \in D''_m$, we have $|f_{uv}(y)| \geq C_m$.

Let us prove that $f_{uv}(y) \geq C_m$. Assume the converse. Then we have $f_{uv}(y) \leq -C_m$, and therefore,

$$f_{\ell v}(y) = f_{\ell u}(y) + f_{uv}(y) \leq (p-1)C_{m-1} - C_m < 0.$$

This contradicts the fact that the column $\ell$ has the largest number of particles.

Suppose that $\tilde{y}$ is obtained from $y$ by adding $pC_{m-1}$ particles to the column $u$. It is easy to see that for $\tilde{y}$ the column $u$ has at least one more particle than any other column.

To complete the proof of the lemma, write

$$\begin{aligned}
& E(f(\zeta_{t+k_m}) - f(\zeta_t) | \zeta_t = y) \\
&= E(f(\zeta_{t+k_m-1}) - f(\zeta_t) | \zeta_t = y) \\
&\quad + \sum_{z \in Y} E(f(\zeta_{t+k_m}) - f(\zeta_{t+k_m-1}) | \zeta_{t+k_m-1} = z) Q_{yz}^{(k_m-1)}.
\end{aligned}$$

By Lemma 2.1, the first term of the right-hand side is bounded above by $p(k_m - 1)$. For the second term write

$$\begin{aligned}
& \sum_{z \in Y} E(f(\zeta_{t+k_m}) - f(\zeta_{t+k_m-1}) | \zeta_{t+k_m-1} = z) Q_{yz}^{(k_m-1)} \\
&= \sum_{z \in Y} \sum_{\substack{\{i,j\} \in E, \\ \{i,j\} \neq \{u,v\}}} E(f_{ij}^2(\zeta_{t+k_m}) - f_{ij}^2(\zeta_{t+k_m-1}) | \zeta_{t+k_m-1} = z) Q_{yz}^{(k_m-1)} \\
&\quad + \sum_{z \in Y} E(f_{uv}^2(\zeta_{t+k_m}) - f_{uv}^2(\zeta_{t+k_m-1}) | \zeta_{t+k_m-1} = z) Q_{yz}^{(k_m-1)},
\end{aligned}$$

where, by Lemma 2.1, the first term of the right-hand side is bounded above by $p - 1$, while the second term can be written as

$$\sum_{z \in Y, z \neq \tilde{y}} E(f_{uv}^2(\zeta_{t+k_m}) - f_{uv}^2(\zeta_{t+k_m-1}) | \zeta_{t+k_m-1} = z) Q_{yz}^{(k_m-1)}$$

$$+ \mathrm{E}(f_{uv}^2(\zeta_{t+k_m}) - f_{uv}^2(\zeta_{t+k_m-1})|\zeta_{t+k_m-1} = \tilde{y})Q_{y\tilde{y}}^{(k_m-1)},$$

which, by Lemma 2.1 and conditions (iii) and (iv) of Theorem 1.3, is less than or equal to

$$\sum_{z \in Y, z \neq \tilde{y}} Q_{yz}^{(k_m-1)} + (1 - 2C_m M)Q_{y\tilde{y}}^{(k_m-1)} \leq 1 - 2C_m M \delta^{k_m-1}.$$

The last inequality follows from condition (ii) of Theorem 1.3 and the fact that it is possible for $\zeta_t$ to reach $\tilde{y}$ from $y$ in $k_m - 1$ steps.

Putting these bounds together, we get

$$\mathrm{E}(f(\zeta_{t+k_m}) - f(\zeta_t)|\zeta_t = y) \leq p(k_m - 1) + (p - 1) + 1 - 2C_m M \delta^{k_m-1},$$

which is less than or equal to $-1$ by our choice of the sequences $k_m$ and $C_m$.
□

We can now complete the proof of Theorem 1.3: Let $A = D_0$, which is finite. Then define the function $k$ as follows: $k(y) = 1$ if $y \in D_0 \cup D_1$ and $k(y) = k_m$ if

$$y \in D_m \setminus (D_0 \cup D_1 \cup \cdots \cup D_{m-1}), \qquad m = 2, \ldots, p.$$

Inequality (2.2) follows from the fact that the function $k$ is bounded and the transition matrix $Q$ is such that, for any $x$, $Q(x, y) = 0$ except for finitely many values of $y$. Since inequality (2.1) follows from Lemma 2.3, we see that the hypotheses of Theorem 2.1 are fulfilled.

REMARK. Theorem 1.3 is still valid if we replace condition (iv) by the following condition:

(iv′) Suppose that $y \in Y$, $i \in V$ and for all $\ell$ such that $\{i, \ell\} \in E$, we have $f_{i\ell}(y) < 0$. Then

$$Q(y, y + e_i) \geq Q(y, y + e_\ell) + M$$

for any $\ell$ as above.

The proof of the theorem requires the following minor modifications:

(a) In the proof of Lemma 2.2 the column $u$ is in the set of the columns that have the minimal number of particles.

(b) In the proof of Lemma 2.3 let $\ell$ be one of the columns having the smallest number of particles. Then, $f_{uv}(y) \leq -C_m$. Suppose that $\tilde{y}$ is obtained from $y$ by adding at least $pC_{m-1}$ particles to each column that is a neighbor of $u$ and belongs to the set $O$. Then for $\tilde{y}$ the column $u$ has at least one particle less than any neighbor column. Also, we take an increasing sequence of integers $C_m$ such that

$$C_1 \geq \frac{1+p}{2M},$$

$$C_m \geq \max\left\{pC_{m-1}, \frac{1+p+p^3 C_{m-1}}{M\delta p^2 C_{m-1}}\right\}, \qquad m = 2, \ldots, p,$$

and let

$$k_m = \begin{cases} 1, & \text{if } m = 1, \\ 1 + p^2 C_{m-1}, & \text{if } m = 2, \ldots, p. \end{cases}$$

**3. Exponential decay of the invariant distribution.** We start this section introducing some notation.

For $0 \leq s \leq t$ and $1 \leq j \leq n$, let

$$\Delta_j^{s,t}(X^n) = X_j^n(t) - X_j^n(s)$$

and for $1 \leq j \leq n-1$, let

$$D_j(X^n(t)) = \sup_{1 \leq i \leq j} \Delta_i X^n(t).$$

Let $P_X$ be the probability associated to this process with initial condition $X$. If this initial condition is $X^n(0) = (0, \ldots, 0)$, we will write $P_0$.

Throughout this section $T_a, T_b'$ and $T_c''$ will be independent Poisson r.v. with parameters $a$, $b$ and $c$, respectively, and $C_i$ and $K_i$ will be strictly positive constants.

Coupling techniques will be often used in this section. These techniques allow us to compare two different processes or two versions of the same process starting from different initial configurations and we will assume that the reader is familiar with them.

3.1. *Construction of the processes, coupling.* The inductive argument used in this section requires a simultaneous construction on the same probability space of the processes for all values of $n$ and for the two types of boundary conditions. This is done as follows:

Let $(N_{0,j}, N_{1,j}, N_{2,j}; j \in \mathbb{N})$ be a collection of independent Poisson processes whose parameters are $\beta_0$ for $N_{0,j}$, $\beta_1 - \beta_0$ for $N_{1,j}$ and $\beta_2 - \beta_1$ for $N_{2,j}$. We let $X_j^n(\cdot)$ increase by 1 at time $t$ if one of the following conditions are satisfied:

(1) $r(X_{j-1}^n(t-), X_j^n(t-), X_{j+1}^n(t-)) = \beta_0$ and $N_{0,j}$ jumps at time $t$,
(2) $r(X_{j-1}^n(t-), X_j^n(t-), X_{j+1}^n(t-)) = \beta_1$ and $N_{0,j} + N_{1,j}$ jumps at time $t$,
(3) $r(X_{j-1}^n(t-), X_j^n(t-), X_{j+1}^n(t-)) = \beta_2$ and $N_{0,j} + N_{1,j} + N_{2,j}$ jumps at time $t$.

The following lemma is an immediate consequence of our simultaneous construction of all processes. It allows us to compare two processes with different and not necessarily deterministic initial conditions.

LEMMA 3.1. *Suppose $X^n(\cdot)$ and $\tilde{X}^n(\cdot)$ are versions of the process with zero boundary conditions constructed as above (with the same Poisson processes):*

(a) *If $\mathrm{P}(X_j^n(0) \geq \tilde{X}_j^n(0)) = 1$ for all $1 \leq j \leq n$, then*
$$\mathrm{P}(X_j^n(t) \geq \tilde{X}_j^n(t), \forall\, t \geq 0) = 1 \qquad \textit{for all } 1 \leq j \leq n.$$

(b) *If $\mathrm{P}(X_j^n(0) - \tilde{X}_j^n(0) = k) = 1$ for all $1 \leq j \leq n$, then*
$$\mathrm{P}(X_j^n(t) - \tilde{X}_j^n(t) = k, \forall\, t \geq 0) = 1 \qquad \textit{for all } 1 \leq j \leq n.$$

REMARK. Suppose that $j < n$, $t > 0$ and let $X \in \{\mathbb{Z}_+\}^n$ be such that $\Delta_j X \geq 0$. Let $X^n(\cdot)$ be the process with zero boundary conditions starting from $X$ and let $X^j(\cdot)$ be the process with zero boundary conditions starting from the restriction of $X$ to the first $j$ coordinates. Then, with the above simultaneous construction of $X^n(\cdot)$ and $X^j(\cdot)$, we have
$$\{\Delta_j X^n(s) \geq 0, \forall\, 0 \leq s \leq t\} \subset \{X_i^n(s) = X_i^j(s), \forall\, 0 \leq s \leq t, i \leq j\}.$$

3.2. *An auxiliary process and some exponential martingales.* We introduce an auxiliary process on $\{\mathbb{Z}_+\}^r \times \{\mathbb{Z}_+\}^r \times \{\mathbb{Z}_+\}^{r+1}$, which we denote
$$(X^r, Z^r, X^{r+1}).$$
The first and third marginal of this process evolve as processes with zero boundary conditions in $\{\mathbb{Z}_+\}^r$ and in $\{\mathbb{Z}_+\}^{r+1}$ respectively and their jumps follow the same collection of Poisson processes $(N_{0,j}, N_{1,j}, N_{2,j}; j \in \mathbb{N})$. The second marginal performs the same jumps as $X^r$ and at any given time $t$, also increases all its coordinates simultaneously by one unit when one of the following two conditions is satisfied:

(i) $X_{r+1}^{r+1}(t-) > X_r^{r+1}(t-)$, $X_{r-1}^{r+1}(t-) \leq X_r^{r+1}(t-)$ and the process $N_{1,r}$ jumps at time $t$.

(ii) $X_{r+1}^{r+1}(t-) > X_r^{r+1}(t-)$, $X_{r-1}^{r+1}(t-) > X_r^{r+1}(t-)$ and the process $N_{2,r}$ jumps at time $t$.

When any of these conditions is satisfied, it may happen that the $r$th coordinate of the $X^r$ process also jumps at time $t$. If that is the case, it is understood that the $Z^r$ process performs both jumps. This means that $Z_r^r$ increases by two units, while all the other coordinates of $Z^r$ increase by one unit.

The following lemma is a consequence from the construction of this process.

LEMMA 3.2. *Suppose the auxiliary process starts with all its $3r+1$ coordinates equal to 0, then*

(3.1) $$\mathrm{P}(X_i^r(t) \leq X_i^{r+1}(t) \leq Z_i^r(t), \forall 1 \leq i \leq r, t \geq 0) = 1$$

*and*

(3.2) $$\mathrm{P}(Z_i^r(t) - X_i^r(t) = Z_1^r(t) - X_1^r(t), \forall 2 \leq i \leq r, t \geq 0) = 1.$$

PROOF. One just needs to check that if the initial condition is such that

$$X_i^r(0) \leq X_i^{r+1}(0) \leq Z_i^r(0) \qquad \forall 1 \leq i \leq r$$

and

$$Z_i^r(0) - X_i^r(0) = Z_1^r(0) - X_1^r(0) \qquad \forall 2 \leq i \leq r,$$

then no jump of the process can break any of these inequalities. This is straightforward except for the inequality involving the $r$th coordinate of $Z^r$ and $X^{r+1}$. However, since the coordinates of $X^{r+1}$ can only make jumps of one unit, we may assume that $X_r^{r+1} = Z_r^r$ and when this happens, conditions (i) and (ii) guarantee that any jump of the $r$th coordinate of $X^{r+1}$ due to its $(r+1)$st coordinate is simultaneous to a jump of all the coordinates of $Z^r$. □

We will later need the following:

LEMMA 3.3. *Suppose the auxiliary process starts with all its $3r+1$ coordinates equal to 0, let*

$$u_s = r(X_{r-1}^{r+1}(s), X_r^{r+1}(s), X_{r+1}^{r+1}(s)) - r(X_{r-1}^{r+1}(s), X_r^{r+1}(s), 0),$$

*and let*

$$v_s = r(X_r^{r+1}(s), X_{r+1}^{r+1}(s), 0),$$

*then, for all $\alpha \geq 0$,*

(3.3) $$(1+\alpha)^{Z_r^r(t) - X_r^r(t)} \exp\left(-\alpha \int_0^t u_s \, ds\right)$$

*and*

(3.4) $$(1+\alpha)^{X_{r+1}^{r+1}(t)} \exp\left(-\alpha \int_0^t v_s \, ds\right)$$

*are martingales.*

PROOF. We will apply Theorem T2 in page 165 of [1] and use as much as possible the notation of that reference. Consider the point process $Z_r^r(t) - X_r^r(t)$ and let $\mathcal{F}_t$ be the $\sigma$-algebra generated by the Poisson processes $N_{i,j}$, $0 \le i \le 2$, $1 \le j \le r+1$. Then, the intensity of this point process is the $\mathcal{F}_t$-predictable process $\lambda_t = \lim_{s \uparrow t} u_s$. Let $\mu_s \equiv 1+\alpha$. Since $\lambda$ and $\mu$ are bounded, condition (2.1) in that reference is satisfied. It then follows from the referred theorem that

$$(3.5) \qquad (1+\alpha)^{Z_r^r(t)-X_r^r(t)} \exp\left(-\alpha \int_0^t u_s\,ds\right)$$

is a local martingale. Since for $0 \le s \le t$, $Z_r^r(s) - X_r^r(s)$ is increasing, positive and bounded by a Poisson random variable of parameter $(\beta_2 - \beta_0)t$, it is also a martingale. A similar argument shows that (3.4) is a martingale too. □

3.3. *Zero boundary conditions.* In this subsection we consider the process with zero boundary conditions. In several parts of our proofs the following observations will play an important role. Suppose $n \in \mathbb{N}$, $0 < i < n$, and $0 \le s < t$, then

$$(3.6) \quad \{\Delta_1 X^n(u) \ge 0, \forall u \in [s,t]\} \subset \{\Delta_1^{s,t}(X^n) = N_{0,1}(t) - N_{0,1}(s)\}$$

and

$$(3.7) \quad \begin{aligned}&\{\Delta_i X^n(u) > 0, \forall u \in [s,t]\} \\ &\quad \subset \{\Delta_{i+1}^{s,t}(X^n) \ge N_{0,i+1}(t) - N_{0,i+1}(s) + N_{1,i+1}(t) - N_{1,i+1}(s)\}.\end{aligned}$$

The strategy of the proof is to proceed by induction in the number of coordinates. The initial lemma provides the result we need when we have two coordinates. The second lemma takes advantage of the boundary effects to show that the first and last coordinates are unlikely to be much higher than there respective neighbors. This provides the initial statement for the induction in index $i$ performed in the proof of Theorem 3.1.

LEMMA 3.4. *There exist $\alpha_2$, $\gamma_2 > 0$, and $0 < d_2 < \beta_1$ such that*

$$P_0(|\Delta_1 X^2(t)| \ge k) \le \exp(-\alpha_2 k) \qquad \forall k \in \mathbb{N}, t \ge 0,$$

*and*

$$P_0(X_2^2(t) \ge d_2 t) \le \exp(-\gamma_2 t) \qquad \forall t \ge 0.$$

PROOF. First note that $|\Delta_1 X^2(t)|$ is a continuous time birth and death process on $\mathbb{Z}_+$, such that:

- $0 \to 1$ at rate $2\beta_0$,
- $n \to n+1$ at rate $\beta_0$ if $n \ge 1$,
- $n \to n-1$ at rate $\beta_1$ if $n \ge 1$.

Since $\beta_0 < \beta_1$, this birth and death process is positive recurrent and its invariant measure has an exponentially decaying tail. It is then easy to couple two versions of this process, one starting from its invariant distribution and the other one from the point mass at 0 in such a way that, with probability 1, the former is at all times above the latter. This proves the first assertion of the lemma. For the second assertion, note that the process $X_1^2(t) + X_2^2(t)$ increases by one unit at a rate which is bounded above by $\beta_0 + \beta_1$. Therefore, it can be coupled with a Poisson process $Z(t)$ of parameter $\beta_0 + \beta_1$ in such a way that

$$P(X_1^2(t) + X_2^2(t) \leq Z(t)) = 1.$$

The second assertion then follows from the first assertion, the inequality

$$2X_2^2(t) \leq X_1^2(t) + X_2^2(t) + |\Delta_1 X^2(t)|$$

and standard large deviations estimates for $Z(t)$. □

LEMMA 3.5. *There exist $C'$ and $\alpha' > 0$ such that, for all $n \in \mathbb{N}$, $k \in \mathbb{N}$, and $t \geq 0$, we have*

(3.8) $$P_0(\Delta_1 X^n(t) \geq k) \leq C' \exp(-\alpha' k)$$

*and*

(3.9) $$P_0(-\Delta_{n-1} X^n(t) \geq k) \leq C' \exp(-\alpha' k).$$

PROOF. Since the process starting from $X \equiv 0$ is invariant under the map

$$(x_1, x_2, \ldots, x_{n-1}, x_n) \to (x_n, x_{n-1}, \ldots, x_2, x_1),$$

under $P_0$ the random variables $\Delta_1 X^n(t)$ and $-\Delta_{n-1} X^n(t)$ have the same distribution. Hence, it suffices to prove (3.8). We assume without loss of generality that $k \geq 4\beta_1$, since this additional condition can be dropped adjusting the constants $C'$ and $\alpha'$. Moreover, we also assume that

(3.10) $$t > \frac{k}{2\beta_1}$$

because when this inequality fails, we have

$$P_0(\Delta_1 X^n(t) \geq k) \leq P((N_{0,1} + N_{1,1})(t) \geq k) \leq P(T_{k/2} \geq k),$$

which decays exponentially in $k$. Let

$$\tau_t = \sup\{s : 0 \leq s \leq t \text{ such that } \Delta_1 X^n(s) = 0\},$$

and let

$$\ell = \min\{q \in \mathbb{N} : (t-q)\beta_1 \leq k/2\}.$$

It follows from this definition and our assumptions on $k$ and $t$ that

$$t - \ell \geq 1,$$

and

(3.11) $$2(t-\ell)\beta_1 \leq k \leq 2(t-\ell+1)\beta_1 \leq 4(t-\ell)\beta_1.$$

We will show that both

$$P_0(\Delta_1 X^n(t) \geq k, \tau_t > \ell)$$

and

$$P_0(\Delta_1 X^n(t) \geq k, \tau_t \leq \ell)$$

decay exponentially in $k$ with constants which depend on $\beta_0$ and $\beta_1$ but not on $t$ or $n$.

For the first of these terms, write

$$\begin{aligned}
P_0(\Delta_1 X^n(t) &\geq k, \tau_t > \ell) \\
&\leq P_0(\Delta_1^{\ell,t}(X^n) \geq k) \\
&\leq P_0((N_{0,1} + N_{1,1})(t) - (N_{0,1} + N_{1,1})(\ell) \geq k) \\
&\leq P(T_{\beta_1(t-\ell)} \geq k) \\
&\leq P(T_{k/2} \geq k),
\end{aligned}$$

which decays exponentially in $k$.

For the second term, let $m \leq \ell$ and note that

$$\{\tau_t \in [m-1, m), \Delta_1 X^n(t) > 0\} \subset \{\Delta_1 X^n(s) > 0, \forall\, s \in [m, t]\}.$$

Hence, we deduce from (3.6) that on $\{\tau_t \in [m-1, m), \Delta_1 X^n(t) > 0\}$ we have

$$\begin{aligned}
\Delta_1^{m-1,t}(X^n) &= \Delta_1^{m-1,m}(X^n) + \Delta_1^{m,t}(X^n) \\
&\leq N_{0,1}(t) - N_{0,1}(m) + (N_{0,1} + N_{1,1})(m) \\
&\quad - (N_{0,1} + N_{1,1})(m-1),
\end{aligned}$$

and from (3.7) that

$$\Delta_2^{m,t}(X^n) \geq (N_{0,2} + N_{1,2})(t) - (N_{0,2} + N_{1,2})(m).$$

Since on $\{\Delta_1 X^n(t) > 0, \tau_t \in [m-1, m)\}$ we also have

$$\Delta_1^{m-1,t}(X^n) > \Delta_2^{m,t}(X^n),$$

we can write

$$\begin{aligned}
\mathrm{P}_0(\Delta_1 X^n(t) &\geq k, \tau_t \in [m-1, m)) \\
&\leq \mathrm{P}_0(\Delta_1 X^n(t) > 0, \tau_t \in [m-1, m)) \\
&\leq \mathrm{P}(N_{0,1}(t) - N_{0,1}(m) + (N_{0,1} + N_{1,1})(m) - (N_{0,1} + N_{1,1})(m-1) \\
&\qquad \geq (N_{0,2} + N_{1,2})(t) - (N_{0,2} + N_{1,2})(m)) \\
&= \mathrm{P}(T_{\beta_0(t-m)+\beta_0+\beta_1} \geq T'_{\beta_1(t-m)}),
\end{aligned}$$

which decays exponentially in $t - m$. Adding this over $1 \leq m \leq \ell$, we obtain that

$$\mathrm{P}_0(\Delta_1 X^n(t) \geq k, \tau_t \leq \ell)$$

decays exponentially in $t - \ell$. The lemma now follows from (3.11). $\square$

We will now state and prove the main result of this section. Theorem 1.2 follows immediately from it.

THEOREM 3.1. *For all $n \geq 2$, there exist $\alpha_n > 0$ and $a_n$ such that*

(3.12) $\quad \mathrm{P}_0(|\Delta_i X^n(t)| \geq k) \leq a_n \exp(-\alpha_n k) \qquad \forall t \geq 0, 1 \leq i \leq n-1, k \in \mathbb{N},$

*and there exist $\gamma_n > 0$, $b_n$, and $0 < d_n < \beta_1$ such that*

(3.13) $\qquad \mathrm{P}_0(X_n^n(t) \geq d_n t) \leq b_n \exp(-\gamma_n t) \qquad \forall t \geq 0.$

The second inequality of the conclusion in this theorem says that the coordinates next to the boundary grow at a speed which is strictly smaller than $\beta_1$. This will play an important role in the inductive step: In time intervals on which coordinate $r$ remains higher than coordinate $r+1$, the first $r$ coordinates behave as a process with $r$ coordinates and zero boundary conditions. Hence, coordinate $r$ grows at a rate which is strictly smaller than $\beta_1$, while coordinate $r+1$ grows at least at that rate. This implies that coordinate $r$ is unlikely to remain higher than coordinate $r+1$ for a long time.

PROOF OF THEOREM 3.1. We proceed by induction on $n$. For $n = 2$, (3.12) and (3.13) follow from Lemma 3.4. Suppose that (3.12) and (3.13) hold for $n \in \{1, \ldots, r\}$.

We start proving the following statement: For $1 \leq i \leq r$, there exists $\alpha'_i > 0$ and $c'_i$ such that

(3.14) $\qquad \mathrm{P}_0(\Delta_i X^{r+1}(t) \geq k) \leq c'_i \exp(-\alpha'_i k) \qquad \forall t \geq 0, k \in \mathbb{N}.$

The proof of (3.14) is done by induction on $i$. For $i = 1$, (3.14) holds by Lemma 3.5. We now suppose that (3.14) holds for $i \in \{1, \ldots, j-1\}$, where $1 \leq j-1 \leq r-1$. Note that, for these values of $j$, by the inductive hypothesis on $n$, (3.13) holds for $n=j$. Let $t \geq 1$, let

$$\tau_t = \sup\{s : 0 \leq s \leq t \text{ such that } \Delta_j X^{r+1}(s) = 0\},$$

$\ell \in \mathbb{N} \cap [0,t]$, $k \in \mathbb{N}$, and let $L$ be large enough to satisfy

$$d_j + (r+1)/L < \beta_1$$

[here $d_j$ is the constant in (3.13) with $n=j$]. Then,

$$\begin{aligned}
&\mathrm{P}_0(\Delta_j X^{r+1}(t) > 0, \tau_t \in [\ell-1, \ell)) \\
&\quad \leq \mathrm{P}_0(\Delta_j X^{r+1}(s) > 0, \forall s \in [\ell, t], D_j(X^{r+1}(\ell)) \leq (t-\ell)/L) \\
&\quad + \mathrm{P}_0(\Delta_j X^{r+1}(\ell) > (t-\ell)/L, \tau_t \in [\ell-1, \ell)) \\
&\quad + \mathrm{P}_0(D_{j-1}(X^{r+1}(\ell)) > (t-\ell)/L).
\end{aligned} \qquad (3.15)$$

We will now show that each of the three terms of the right-hand side above is bounded above by expressions of the form $C\exp(-K(t-\ell))$.

For the third term, this follows by the inductive hypothesis on $i$.

For the the second term, note that

$$\{\Delta_j X^{r+1}(\ell) > (t-\ell)/L, \tau_t \in [\ell-1, \ell)\} \subset \{\Delta_j^{\ell-1,\ell}(X^{r+1}) \geq (t-\ell)/L\}.$$

Therefore,

$$\begin{aligned}
&\mathrm{P}_0(\Delta_j X^{r+1}(\ell) > (t-\ell)/L, \tau_t \in [\ell-1, \ell)) \\
&\quad \leq \mathrm{P}((N_{0,j} + N_{1,j} + N_{2,j})(\ell) - (N_{0,j} + N_{1,j} + N_{2,j})(\ell-1) > (t-\ell)/L) \\
&\quad = \mathrm{P}(T_{\beta_2} > (t-\ell)/L),
\end{aligned}$$

which decays exponentially in $t-\ell$.

For the first term, write

$$\begin{aligned}
&\{\Delta_j X^{r+1}(s) > 0, \forall s \in [\ell, t], D_j(X^{r+1}(\ell)) \leq (t-\ell)/L\} \\
&\quad \subset \{\Delta_j X^{r+1}(t) > 0, D_j(X^{r+1}(\ell)) \leq (t-\ell)/L\} \\
&\quad \subset \{\Delta_j X^{r+1}(t) > 0, \Delta_j X^{r+1}(\ell) \leq (t-\ell)/L\} \\
&\quad \subset \{\Delta_j^{\ell,t}(X^{r+1}) > \Delta_{j+1}^{\ell,t}(X^{r+1}) - (t-\ell)/L\} \\
&\quad \subset \{\Delta_j^{\ell,t}(X^{r+1}) > (d_j + r/L)(t-\ell)\} \\
&\quad \qquad \cup \{\Delta_{j+1}^{\ell,t}(X^{r+1}) < (d_j + (r+1)/L)(t-\ell)\}.
\end{aligned}$$

Hence,

$$\begin{aligned}
&\{\Delta_j X^{r+1}(s) > 0, \forall s \in [\ell, t], D_j(X^{r+1}(\ell)) \leq (t-\ell)/L\} \\
&\quad \subset \{D_j(X^{r+1}(\ell)) \leq (t-\ell)/L\} \\
&\quad \cap (\{\Delta_j^{\ell,t}(X^{r+1}) > (d_j + r/L)(t-\ell), \Delta_j X^{r+1}(s) > 0, \forall s \in [\ell, t]\}
\end{aligned}$$

$$\cup \{\Delta_{j+1}^{\ell,t}(X^{r+1}) < (d_j + (r+1)/L)(t-\ell), \Delta_j X^{r+1}(s) > 0,$$
$$\forall s \in [\ell, t]\}).$$

Therefore, applying the Markov property at time $\ell$ and letting

$$A_{m,j} = \{X \in \{\mathbb{Z}_+\}^m : D_j(X) \leq (t-\ell)/L\},$$

where $j < m \in \mathbb{N}$, we get

$$\begin{aligned}
\mathrm{P}_0(\Delta_j X^{r+1}(s) > 0, \forall s \in [\ell, t], D_j(X^{r+1}(\ell)) \leq (t-\ell)/L) \\
\leq \sup_{X \in A_{r+1,j}} \mathrm{P}_X(\Delta_j^{0,t-\ell}(X^{r+1}) > (d_j + r/L)(t-\ell), \\
\Delta_j X^{r+1}(s) > 0, \forall s \in [0, t-\ell]) \\
+ \sup_{X \in A_{r+1,j}} \mathrm{P}_X(\Delta_{j+1}^{0,t-\ell}(X^{r+1}) < (d_j + (r+1)/L)(t-\ell), \\
\Delta_j X^{r+1}(s) > 0, \forall s \in [0, t-\ell]).
\end{aligned} \quad (3.16)$$

From the remark following Lemma 3.1, we get that the first term of the right-hand side of (3.16) is bounded by

$$\sup_{X \in A_{j,j-1}} \mathrm{P}_X(\Delta_j^{0,t-\ell}(X^j) > (d_j + r/L)(t-\ell)),$$

while, by (3.7), the second term is bounded above by

$$\mathrm{P}((N_{0,j+1} + N_{1,j+1})(t-\ell) \leq (d_j + (r+1)/L)(t-\ell)).$$

For $X \in A_{j,j-1}$, let $X'$ be the element of $\{\mathbb{Z}_+\}^j$ whose coordinates are all equal to $\max\{X_1, \ldots, X_j\}$. Note that the value of these coordinates is bounded above by $X_j + (j-1)(t-\ell)/L$. Therefore, by parts (a) and (b) of Lemma 3.1, the first term is bounded above by

$$\begin{aligned}
\mathrm{P}_{X'}(\Delta_j^{0,t-\ell}(X^j) &> (d_j + r/L - (j-1)/L)(t-\ell)) \\
&= \mathrm{P}_0(\Delta_j^{0,t-\ell}(X^j) > (d_j + r/L - (j-1)/L)(t-\ell)) \\
&\leq \mathrm{P}_0(\Delta_j^{0,t-\ell}(X^j) > d_j(t-\ell)) \\
&= \mathrm{P}_0(X_j^j(t-\ell) > d_j(t-\ell)) \\
&\leq b_j \exp(-\gamma_j(t-\ell)),
\end{aligned}$$

where the last inequality follows from the inductive hypothesis on $n$ (recall that $j \leq r$). Since the Poisson process $N_{0,j+1} + N_{1,j+1}$ has parameter $\beta_1 > d_j + (r+1)/L$, we get that the second term of the right-hand side of (3.16) also decays exponentially in $t - \ell$. This completes the proof that the first term of the right-hand side of (3.15) decays exponentially.

Since the three terms of the right-hand side of (3.15) decay exponentially in $t - \ell$, we have proved that there exist constants $C_1$ and $K_1 > 0$ such that

$$P_0(\Delta_j X^{r+1}(t) > 0, \tau_t \in [\ell - 1, \ell)) \leq C_1 \exp(-K_1(t - \ell))$$

holds for all $t \geq 1$, $\ell \leq t$, and $\ell \in \mathbb{N}$. Therefore, substituting $\ell'$ for $\ell$ and summing on $1 \leq \ell' \leq \ell$, we get

(3.17) $$\begin{aligned} P_0(\Delta_j X^{r+1}(t) > 0, \tau_t \leq \ell) \\ \leq C_2 \exp(-K_2(t - \ell)) \end{aligned} \qquad \forall t \geq 1, \ell \leq t, \ell \in \mathbb{N}.$$

We will now complete the inductive step for (3.14). We may assume that

(3.18) $$t \geq \frac{k}{2\beta_2},$$

since otherwise

$$\begin{aligned} P_0(\Delta_j X^{r+1}(t) \geq k) \\ \leq P((N_{0,j} + N_{1,j} + N_{2,j})(t) \geq k) \\ \leq P(T_{k/2} \geq k), \end{aligned}$$

which decays exponentially in $k$ with a constant which does not depend on $t$. And arguing as in the proof of Lemma 3.5, we may also assume that $k > 4\beta_2$. Let

$$\ell = \min\{q \in \mathbb{N} : (t - q)\beta_2 \leq k/2\}.$$

Then, it follows from our assumptions on $k$ and $t$ that

(3.19) $$t - \ell \geq 1$$

and

(3.20) $$k \leq 2(t - \ell + 1)\beta_2 \leq 4(t - \ell)\beta_2.$$

Now write

$$\begin{aligned} P_0(\Delta_j X^{r+1}(t) \geq k) \\ \leq P_0(\Delta_j X^{r+1}(t) > 0, \tau_t \leq \ell) \\ + P_0(\Delta_j X^{r+1}(t) \geq k, \tau_t > \ell). \end{aligned}$$

By (3.17), the first term is bounded above by

$$C_2 \exp(-K_2(t - \ell)) \leq C_2 \exp\left(-\frac{K_2}{4\beta_2}k\right),$$

and the second term is bounded above by

$$P_0(\Delta_j^{\ell,t}(X^{r+1}) \geq k) \leq P(T_{\beta_2(t-\ell)} \geq k) \leq P(T_{k/2} \geq k).$$

Therefore, both terms decay exponentially in $k$ with constants which do not depend on $t$. Thus, (3.14) holds for $i = j$. Hence, by induction, it holds for all $i$, and by symmetry, it follows that (3.12) holds for $n = r + 1$. It remains to prove that (3.13) also holds for that $n$.

Let $u_s$ and $v_s$ be as in Lemma 3.3. Then, for all $s \geq 0$, we have

$$0 \leq u_s \leq \max\{\beta_2 - \beta_1, \beta_1 - \beta_0\} \leq \beta_2 - \beta_0, \qquad 0 \leq v_s \leq \beta_1,$$

and that $u_s > 0$ implies $X_{r+1}^{r+1}(s) > X_r^{r+1}(s)$, which in turn implies $v_s = \beta_0$. Hence, if $0 \leq \delta \leq 1$ and $t \geq 0$, then either

$$\int_0^t u_s \, ds \leq (\beta_2 - \beta_0)\delta t$$

or

$$\int_0^t v_s \, ds \leq [\beta_0 \delta + \beta_1(1 - \delta)]t.$$

Let

$$\gamma \in (0, \beta_1 - d_r),$$

then pick $\delta \in (0, 1)$ and $\eta > 0$ such that

$$\gamma - \delta(\beta_2 - \beta_0) > 0, \qquad \eta + \gamma < \beta_1 - d_r$$

and

$$\eta < \delta(\beta_1 - \beta_0).$$

Let $(X^r, Z^r, X^{r+1})$ be the auxiliary process defined in Section 3.2 and let $P_0$ be the probability associated to this process when it starts with all its coordinates equal to 0. Then, write

$$P_0(Z_r^r(t) - X_r^r(t) \geq \gamma t, X_{r+1}^{r+1}(t) \geq (\beta_1 - \eta)t)$$

$$\leq P_0\left(Z_r^r(t) - X_r^r(t) \geq \gamma t, \int_0^t u_s \, ds \leq (\beta_2 - \beta_0)\delta t\right)$$

$$+ P_0\left(X_{r+1}^{r+1}(t) \geq (\beta_1 - \eta)t, \int_0^t v_s \, ds \leq \beta_0\delta + \beta_1(1 - \delta)t\right)$$

$$\leq P_0\left(Z_r^r(t) - X_r^r(t) - \int_0^t u_s \, ds \geq \gamma - \delta(\beta_2 - \beta_0)t\right)$$

$$+ P_0\left(X_{r+1}^{r+1}(t) - \int_0^t v_s \, ds \geq [\delta(\beta_1 - \beta_0) - \eta]t\right).$$

Let

$$c = \min\{\gamma - \delta(\beta_2 - \beta_0), \delta(\beta_1 - \beta_0) - \eta\},$$

then $c > 0$ and

(3.21)
$$\begin{aligned}
&\mathrm{P}_0(Z_r^r(t) - X_r^r(t) \geq \gamma t, X_{r+1}^{r+1}(t) \geq (\beta_1 - \eta)t) \\
&\quad \leq \mathrm{P}_0\Big(Z_r^r(t) - X_r^r(t) - \int_0^t u_s\, ds \geq ct\Big) \\
&\quad\quad + \mathrm{P}_0\Big(X_{r+1}^{r+1}(t) - \int_0^t v_s\, ds \geq ct\Big).
\end{aligned}$$

Let $\rho > 0$ be such that
$$\varepsilon =: (1-\rho)c - \rho\beta_2 > 0$$
and let $\alpha > 0$ be such that
$$\alpha(1-\rho) \leq \ln(1+\alpha).$$

From Lemma 3.3, we have
$$\mathrm{E}\Big(\exp\Big[\ln(1+\alpha)(Z_r^r(t) - X_r^r(t)) - \alpha \int_0^t u_s\, ds\Big]\Big) = 1.$$

Therefore,
$$\mathrm{E}\Big(\exp\Big[\alpha\Big((1-\rho)(Z_r^r(t) - X_r^r(t)) - \alpha \int_0^t u_s\, ds\Big)\Big]\Big) \leq 1,$$

which, by Chebyshev's inequality, implies that
$$\mathrm{P}\Big((1-\rho)(Z_r^r(t) - X_r^r(t)) - \int_0^t u_s\, ds \geq a\Big) \leq \exp(-\alpha a) \qquad \forall a > 0.$$

Hence,

(3.22)
$$\begin{aligned}
&\mathrm{P}_0\Big((Z_r^r(t) - X_r^r(t)) - \int_0^t u_s\, ds \geq ct\Big) \\
&\quad = \mathrm{P}_0\Big((1-\rho)(Z_r^r(t) - X_r^r(t)) - \int_0^t u_s\, ds \geq (1-\rho)ct - \rho\int_0^t u_s\, ds\Big) \\
&\quad \leq \mathrm{P}_0\Big((1-\rho)(Z_r^r(t) - X_r^r(t)) - \int_0^t u_s\, ds \geq ((1-\rho)c - \rho\beta_2)t\Big) \\
&\quad \leq \exp(-\alpha\varepsilon t).
\end{aligned}$$

Using the second martingale provided by Lemma 3.3,
$$(1+\alpha)^{X_{r+1}^{r+1}(t)} \exp\Big(-\alpha \int_0^t v_s\, ds\Big)$$
and proceeding as above, we get

(3.23) $\quad \mathrm{P}_0\Big(X_{r+1}^{r+1}(t) - \int_0^t v_s\, ds \geq ct\Big) \leq \exp(-\alpha\varepsilon' t),$

where $\varepsilon' = (1-\rho)c - \rho\beta_1 > \varepsilon$. It now follows from (3.21), (3.22) and (3.23) that
$$\mathrm{P}_0(Z_r^r(t) - X_r^r(t) \geq \gamma t, X_{r+1}^{r+1}(t) \geq (\beta_1 - \eta)t) \leq 2\exp(-\alpha\varepsilon t) \qquad \forall t \geq 0.$$

However, by the inductive hypothesis on $n$, we have
$$P_0(X_r^r(t) \geq d_r t) \leq b_n \exp(-\gamma_r t),$$
therefore, there exist $C_6$ and $K_6 > 0$ such that

(3.24) $$\begin{aligned} P_0(Z_r^r(t) \geq (d_r + \gamma)t, X_{r+1}^{r+1}(t) \geq (\beta_1 - \eta)t) \\ \leq C_6 \exp(-K_6 t) \quad \forall t \geq 0. \end{aligned}$$

Since $X_r^{r+1}(t) \leq Z_r^r(t)$, we have
$$P_0(Z_r^r(t) \leq (d_r + \gamma)t, X_{r+1}^{r+1}(t) \geq (\beta_1 - \eta)t)$$
$$\leq P_0(X_r^{r+1}(t) \leq (d_r + \gamma)t, X_{r+1}^{r+1}(t) \geq (\beta_1 - \eta)t)$$
$$\leq P_0(X_{r+1}^{r+1}(t) - X_r^{r+1}(t) \geq (\beta_1 - \eta - \gamma - d_r)t).$$

Since $\beta_1 - \eta - \gamma - d_r > 0$, it follows from Lemma 3.5 that there exist constants $C_7$ and $K_7 > 0$ such that

(3.25) $$\begin{aligned} P_0(Z_r^r(t) \leq (d_r + \gamma)t, X_{r+1}^{r+1}(t) \geq (\beta_1 - \eta)t) \\ \leq C_7 \exp(-K_7 t) \quad \forall t \geq 0. \end{aligned}$$

By (3.24) and (3.25), we get $C_8$ and $K_8 > 0$ such that
$$P_0(X_{r+1}^{r+1}(t) \geq (\beta_1 - \eta)t) \leq C_8 \exp(-K_8 t) \quad \forall t \geq 0.$$

Hence, (3.13) holds for $n = r + 1$ and the induction in $n$ is complete. □

3.4. *Periodic boundary conditions.* We denote by
$$Y^n(t) = (Y_1^n(t), \ldots, Y_n^n(t))$$
the process with periodic boundary conditions on $\{\mathbb{Z}_+\}^n$.

For this process, we define for $1 \leq i \leq n$
$$\Delta_i Y^n(t) = Y_i^n(t) - Y_{i+1}^n(t),$$
where, by convention, $Y_{n+1}^n(t) = Y_1^n(t)$. Note that the process
$$(\Delta_1 Y^n(t), \ldots, \Delta_n Y^n(t))$$
is Markovian too and $\sum_{i=1}^n \Delta_i Y^n(t) = 0$.

For $0 \leq s \leq t$ and $1 \leq j \leq n$, let
$$\Delta_j^{s,t}(Y^n) = Y_j^n(t) - Y_j^n(s).$$

Since $\sum_{i=1}^n \Delta_i Y^n(t) = 0$ and the semigroup of the process commutes with the map
$$(y_1, y_2, \ldots, y_n) \to (y_2, \ldots, y_n, y_1),$$
Theorem 1.1 follows from the following:

THEOREM 3.2. *For all $1 \leq i \leq n-1$, there exist positive $\bar{C}_i$ and $\bar{K}_i > 0$ such that*

$$P_0(\min\{-\Delta_n Y^n(t), \Delta_i Y^n(t)\} \geq k) \leq \bar{C}_i \exp(-\bar{K}_i k) \qquad \forall t \geq 0, k \in \mathbb{N}.$$

The idea of the proof is that as long as $-\Delta_n$ and $\Delta_i$ are positive, coordinates $1,\ldots,i$ behave as a process with zero boundary conditions. By Theorem 3.1, they grow at a rate which is strictly smaller than $\beta_1$. But coordinates $i+1$ and $n$ grow at least at that rate. Therefore, $-\Delta_n$ and $\Delta_i$ are unlikely to remain positive too long.

PROOF OF THEOREM 3.2. As in the proof of Theorem 3.1, it suffices to prove the result under the additional assumptions $k \geq 4\beta_2$ and $t > \frac{k}{2\beta_2}$. We proceed by induction on $i$.

*First step.* $i = 1$. Let

$$\tau_t = \sup\{s \in [0,t] : \Delta_n Y^n(s) = 0 \text{ or } \Delta_1 Y^n(s) = 0\},$$

and let

$$\ell = \min\{q \in \mathbb{N} : (t-q)\beta_2 \leq k/2\}.$$

Then

$$\begin{aligned}
P_0(\min\{-\Delta_n Y^n(t), \Delta_1 Y^n(t)\} \geq k, \tau_t \geq \ell) \\
\leq P((N_{0,1} + N_{1,1} + N_{2,1})(t) - (N_{0,1} + N_{1,1} + N_{2,1})(\ell) \geq k) \\
= P(T_{\beta_2(t-\ell)} \geq k) \leq P(T_{k/2} \geq k) \\
\leq C_5 \exp(-K_5 k).
\end{aligned} \tag{3.26}$$

For $m \leq \ell$, write

$$P_0(\min\{-\Delta_n Y^n(t), \Delta_1 Y^n(t)\} \geq k, \tau_t \in [m-1,m))$$
$$\leq P_0(\min\{-\Delta_n Y^n(s), \Delta_1 Y^n(s)\} > 0, \forall s \in [m,t], \tau_t \in [m-1,m)).$$

But on the set

$$\{\min\{-\Delta_n Y^n(s), \Delta_1 Y^n(s)\} > 0, \forall s \in [m,t], \tau_t \in [m-1,m)\}$$

we have

$$N_{0,1}(t) - N_{0,1}(m-1) \geq \Delta_1^{\tau_t,t}(Y^n) \geq \min\{\Delta_2^{\tau_t,t}(Y^n), \Delta_n^{\tau_t,t}(Y^n)\}$$
$$\geq \min\{\Delta_2^{m,t}(Y^n), \Delta_n^{m,t}(Y^n)\}.$$

Since on that set we also have

$$\Delta_2^{\tau_t,t}(Y^n) \geq (N_{0,2} + N_{1,2})(t) - (N_{0,2} + N_{1,2})(m)$$

and
$$\Delta_n^{\tau_t,t}(Y^n) \geq (N_{0,n} + N_{1,n})(t) - (N_{0,n} + N_{1,n})(m),$$
we get
$$\begin{aligned}
P_0(\min\{-\Delta_n Y^n(t), \Delta_1 Y^n(t)\} \geq k, \tau_t \in [m-1, m)) \\
\leq P(N_{0,1}(t) - N_{0,1}(m-1) \\
\geq [(N_{0,2} + N_{1,2})(t) - (N_{0,2} + N_{1,2})(m)] \\
\wedge [(N_{0,n} + N_{1,n})(t) - (N_{0,n} + N_{1,n})(m)]) \\
= P(T_{\beta_0(t-m+1)} > T'_{\beta_1(t-m)} \wedge T''_{\beta_1(t-m)}) \\
\leq 2 P(T_{\beta_0(t-m+1)} > T'_{\beta_1(t-m)}).
\end{aligned}$$
As $\beta_0 < \beta_1$, this decays exponentially in $t-m$. Thus, for some constants $C_6$ and $K_6 > 0$, we have
$$\begin{aligned}
P_0(\min\{-\Delta_n Y^n(t), \Delta_1 Y^n(t)\} \geq k, \tau_t \in [m-1, m)) \\
\leq C_6 \exp(-K_6(t-m)).
\end{aligned}$$
Summing over $1 \leq m \leq \ell$, we get
$$P_0(\min\{-\Delta_n Y^n(t), \Delta_1 Y^n(t)\} \geq k, \tau_t < \ell) \leq C_7 \exp(-K_7(t-\ell)).$$
This, (3.20) and (3.26) imply the result for $i=1$.

*Second step.* Suppose the result holds for $1 \leq i \leq r-1 < n-1$. Now let
$$\tau_t = \sup\{s \in [0,t] : \Delta_n Y^n(t) = 0 \text{ or } \Delta_r Y^n(t) = 0\}$$
and as in the first step, let
$$\ell = \min\{q \in \mathbb{N} : (t-q)\beta_2 \leq k/2\}.$$
Then
$$\begin{aligned}
P_0(\min\{-\Delta_n Y^n(t), \Delta_r Y^n(t)\} \geq k, \tau_t \geq \ell) \\
\leq P((N_{0,1} + N_{1,1} + N_{2,1})(t) - (N_{0,1} + N_{1,1} + N_{2,1})(\ell) \geq k) \\
+ P((N_{0,r} + N_{1,r} + N_{2,r})(t) - (N_{0,r} + N_{1,r} + N_{2,r})(\ell) \geq k) \\
\leq 2 P(T_{\beta_2(t-\ell)} \geq k) \leq 2 P(T_{k/2} \geq k) \\
\leq C_8 \exp(-K_8 k).
\end{aligned} \quad (3.27)$$
Let $m \in \mathbb{N}$ be less than or equal to $\ell$. We wish to show that
$$P_0(\min\{-\Delta_n Y^n(t), \Delta_r Y^n(t)\} \geq k, \tau_t \in [m-1, m))$$
decays exponentially in $t-m$. To do so, first note that if
$$\min\{-\Delta_n Y^n(t), \Delta_r Y^n(t)\} \geq k,$$

then either $-\Delta_n Y^n(\tau_t) = 1$ or $\Delta_r Y^n(\tau_t) = 1$. Therefore,

$$\begin{aligned}
&\mathrm{P}_0(\min\{-\Delta_n Y^n(t), \Delta_r Y^n(t)\} \geq k, \tau_t \in [m-1, m)) \\
&\quad \leq \mathrm{P}_0(\min\{-\Delta_n Y^n(t), \Delta_r Y^n(t)\} \geq k, \tau_t \in [m-1, m), -\Delta_n Y^n(\tau_t) = 1) \\
&\quad\quad + \mathrm{P}_0(\min\{-\Delta_n Y^n(t), \Delta_r Y^n(t)\} \geq k, \tau_t \in [m-1, m), \Delta_r Y^n(\tau_t) = 1) \\
&\quad = 2\,\mathrm{P}_0(\min\{-\Delta_n Y^n(t), \Delta_r Y^n(t)\} \geq k, \tau_t \in [m-1, m), \Delta_r Y^n(\tau_t) = 1),
\end{aligned}$$

where the equality follows from the facts that the process is invariant under the map

$$(y_1, \ldots, y_n) \to (y_r, y_{r-1}, \ldots, y_1, y_n, y_{n-1}, \ldots, y_{r+1})$$

and the initial configuration is a fixed point of that map.

Let $\alpha > 0$ be such that $2\alpha + d_r < \beta_1$ ($d_r$ is the same constant as in Theorem 3.1) and let

$$\begin{aligned}
E_0 &= \{\min\{-\Delta_n Y^n(t), \Delta_r Y^n(t)\} \geq k, \tau_t \in [m-1, m), \Delta_r Y^n(\tau_t) = 1\}, \\
E_1 &= \{\max\{Y_1^n(\tau_t), \ldots, Y_r^n(\tau_t)\} \geq \max\{Y_1^n(\tau_t), Y_r^n(\tau_t)\} + \alpha(t-m)\}, \\
E_2 &= \{\max\{Y_1^n(\tau_t), \ldots, Y_r^n(\tau_t)\} \leq Y_r^n(\tau_t) + \alpha(t-m)\}, \\
E_3 &= \{\max\{Y_1^n(\tau_t), \ldots, Y_r^n(\tau_t)\} \leq Y_1^n(\tau_t) + \alpha(t-m), \\
&\qquad\qquad -\Delta_n Y^n(\tau_t) < \alpha(t-m)\} \\
E_4 &= \{Y_r^n(\tau_t) + \alpha(t-m) \leq \max\{Y_1^n(\tau_t), \ldots, Y_r^n(\tau_t)\} \leq Y_1^n(\tau_t) + \alpha(t-m), \\
&\qquad\qquad -\Delta_n Y^n(\tau_t) \geq \alpha(t-m)\}.
\end{aligned}$$

Since at least one of the last four events must occur, we have

(3.28) $\mathrm{P}_0(E_0) \leq \mathrm{P}_0(E_0 \cap E_1) + \mathrm{P}_0(E_0 \cap E_2) + \mathrm{P}_0(E_0 \cap E_3) + \mathrm{P}_0(E_0 \cap E_4).$

We will now show that the four terms of the right-hand side of (3.28) decay exponentially in $t - m$.

*First term.* This term is bounded above by

(3.29)
$$\begin{aligned}
&\mathrm{P}_0(\max\{Y_1^n(\tau_t), \ldots, Y_r^n(\tau_t)\} \geq \max\{Y_1^n(\tau_t), Y_r^n(\tau_t)\} + \alpha(t-m), \\
&\qquad\qquad\qquad\qquad \tau_t \in [m-1, m)) \\
&\quad \leq \mathrm{P}_0(\max\{Y_1^n(m), \ldots, Y_r^n(m)\} \geq \max\{Y_1^n(m), Y_r^n(m)\} + \alpha(t-m)/2) \\
&\quad\quad + \mathrm{P}(\Delta_1^{m-1,m}(Y^n) \geq \alpha(t-m)/2) + \mathrm{P}(\Delta_r^{m-1,m}(Y^n) \geq \alpha(t-m)/2).
\end{aligned}$$

But

$$\max\{Y_1^n(m), \ldots, Y_r^n(m)\} \geq \max\{Y_1^n(m), Y_r^n(m)\} + \alpha(t-m)/2$$

implies that there exist $1 \leq i_1 < i_2 \leq r-1$ such that

$$-\Delta_{i_1} Y^n(m) \geq \frac{\alpha}{2r}(t-m)$$

and

$$\Delta_{i_2} Y^n(m) \geq \frac{\alpha}{2r}(t-m).$$

Hence, the first term of the right-hand side of (3.29) decays exponentially in $t-m$ by the inductive hypothesis and the invariance of the process under the map

$$(y_1, \ldots, y_n) \to (y_{i_1+1}, \ldots, y_n, y_1, \ldots, y_{i_1}).$$

Since the second and third terms trivially share this property, the same happens to the left-hand side of (3.29).

*Second term.* On the set $A = E_0 \cap E_2$, the coordinates $Y_1^n, \ldots, Y_r^n$ jump in the time interval $(\tau_t, t]$ at the same rates as the coordinates of the process with zero boundary conditions $(X_1^r, \ldots, X_r^r)$. Therefore, if we start this last process at time $m-1$ with all its coordinates equal to

$$\max\{Y_1^n(\tau_t), \ldots, Y_r^n(\tau_t)\},$$

we can couple it with the $Y^n$ process in such a way that on the set $A$ we have $Y_i^n(s) \leq X_i^r(s)$ for all $s \in [\tau_t, t]$ and all $1 \leq i \leq r$. Therefore, on $A$ we have

$$X_r^r(t) \geq Y_r^n(t) \geq Y_{r+1}^n(t) + k \geq Y_{r+1}^n(t) + 1$$

and

$$\begin{aligned} X_r^r(m-1) &= \max\{Y_1^n(\tau_t), \ldots, Y_r^n(\tau_t)\} \\ &\leq Y_r^n(\tau_t) + \alpha(t-m) \\ &= Y_{r+1}^n(\tau_t) + 1 + \alpha(t-m) \\ &\leq Y_{r+1}^n(m) + 1 + \alpha(t-m), \end{aligned}$$

which imply

$$\Delta_r^{m-1,t}(X^r) + \alpha(t-m) \geq \Delta_{r+1}^{m,t}(Y^n).$$

Since on $A$

$$\Delta_{r+1}^{m,t}(Y^n) \geq (N_{0,r+1} + N_{1,r+1})(t) - (N_{0,r+1} + N_{1,r+1})(m),$$

we must also have

$$\Delta_r^{m-1,t}(X^r) + \alpha(t-m) \geq (N_{0,r+1} + N_{1,r+1})(t) - (N_{0,r+1} + N_{1,r+1})(m).$$

Hence,

$$A \subset A_1 \cup A_2,$$

where
$$A_1 = \{\Delta_r^{m-1,t}(X^r) \geq d_r(t-m)\}$$
and
$$A_2 = \{(N_{0,r+1} + N_{1,r+1})(t) - (N_{0,r+1} + N_{1,r+1})(m) \leq (d_r + \alpha)(t-m)\}.$$
Since $\Delta_r^{m-1,t}(X^r)$ is distributed as $X_r^r(t-m+1)$ under $P_0$, the probability of $A_1$ decays exponentially in $t-m$ by Theorem 3.1, and the probability of $A_2$ does the same because $\beta_1 > d_r + 2\alpha > d_r + \alpha$.

*Third term.* Let
$$B = \{\min\{-\Delta_n Y^n(t), \Delta_r Y^n(t)\} \geq k, \tau_t \in [m-1, m)\} \cap E_3.$$
On $B$ the coordinates $Y_1^n, \ldots, Y_r^n$ jump in the time interval $(\tau_t, t]$ at the same rates as the coordinates of the process with zero boundary conditions $(X_1^r, \ldots, X_r^r)$. Therefore, if we start this last process at time $m-1$ with all its coordinates equal to $\max\{Y_1^n(\tau_t), \ldots, Y_r^n(\tau_t)\}$, we can couple it with the $Y^n$ process in such a way that on the set $B$ we have $Y_i^n(s) \leq X_i^r(s)$ for all $s \in [\tau_t, t]$ and all $1 \leq i \leq r$. Therefore, on $B$ we have
$$X_1^r(t) \geq Y_1^n(t) \geq Y_n^n(t)$$
and
$$X_1^r(m-1) = \max\{Y_1^n(\tau_t), \ldots, Y_r^n(\tau_t)\}$$
$$\leq Y_1^n(\tau_t) + \alpha(t-m) \leq Y_n^n(\tau_t) + 2\alpha(t-m)$$
$$\leq Y_n^n(m) + 2\alpha(t-m).$$
Hence, on $B$ we also have
$$\Delta_1^{m-1,t}(X^r) + 2\alpha(t-m) \geq \Delta_n^{m,t}(Y^n).$$
Proceeding as we did for the second term, we get
$$B \subset \{\Delta_1^{m-1,t}(X^r) \geq d_r(t-m)\}$$
$$\cup \{(N_{0,n} + N_{1,n})(t) - (N_{0,n} + N_{1,n})(m) \leq (d_r + 2\alpha)(t-m)\}$$
and both sets have probabilities which decay exponentially in $t-m$.

*Fourth term.* Since
$$\max\{Y_1^n(\tau_t), \ldots, Y_r^n(\tau_t)\} \geq Y_r^n(\tau_t) + \alpha(t-m),$$
there exists $1 \leq i \leq r-1$ such that
$$\Delta_i Y^n(\tau_t) \geq \frac{\alpha}{r}(t-m).$$

Since in this case we also have $-\Delta_n Y^n(\tau_t) \geq \alpha(t-m)$ considering the coordinates $1, \ldots, i$, we can proceed as for the first term to show that the fourth term also decays exponentially in $t - m$.

Hence, there exists constants such that

$$P_0(\min\{-\Delta_n Y^n(t), \Delta_r Y^n(t)\} \geq k, \tau_t \in [m-1, m))$$
$$\leq C_9 \exp(-K_9(t-m)).$$

Adding this on $1 \leq m \leq \ell$, we obtain

$$P_0(\min\{-\Delta_n Y^n(t), \Delta_r Y^n(t)\} \geq k, \tau_t \leq \ell)$$
$$\leq C_{10} \exp(-K_{10}(t-\ell)).$$

This, (3.20) and (3.27) complete the inductive step. $\square$

**Acknowledgment.** We are grateful to the referee for his very careful reading and useful comments.

E. D. ANDJEL
CMI
UNIVERSITÉ DE PROVENCE
39 RUE JOLIOT CURIE
13453 MARSEILLE CEDEX 13
FRANCE
E-MAIL: andjel@cmi.univ-mrs.fr

M. V. MENSHIKOV
DEPARTMENT OF MATHEMATICAL SCIENCES
UNIVERSITY OF DURHAM
SOUTH ROAD, DURHAM DH1 3LE
UK
E-MAIL: Mikhail.Menshikov@durham.ac.uk

V. V. SISKO
DEPARTAMENTO DE ESTATÍSTICA
INSTITUTO DE MATEMÁTICA E ESTATÍSTICA
UNIVERSIDADE DE SÃO PAULO
RUA DO MATÃO 1010
CEP 05508-090, SÃO PAULO, SP
BRAZIL
E-MAIL: valentin@ime.usp.br